\documentclass[]{amsart}
\usepackage[foot]{amsaddr}
\usepackage[a4paper,top=2.5cm,bottom=3cm,inner=3cm,outer=3cm]{geometry}

\usepackage[utf8]{inputenc}

\usepackage{amsmath,amsthm}
\usepackage{amsfonts,amssymb,mathrsfs}

\usepackage{xy}
\usepackage[bookmarks=false,breaklinks=true]{hyperref}
\usepackage{float}
\usepackage{graphicx}
\usepackage[percent]{overpic}
\usepackage{color}

\usepackage{epigraph}
\definecolor{myurlcolor}{rgb}{0,0,0.7}
\usepackage{hyperref}
\hypersetup{colorlinks,
linkcolor=myurlcolor,
citecolor=myurlcolor,
urlcolor=myurlcolor}
\usepackage[mathscr]{euscript}

\input xy
\xyoption{all}

\newcommand{\maps}{\colon}    

\newcommand{\R}{{\mathbb R}}  
\newcommand{\C}{{\mathbb C}}  

\newcommand{\E}{\vec{E}}
\renewcommand{\a}{\vec{a}}
\newcommand{\z}{\vec{z}}
\renewcommand{\v}{\vec{v}}

\newcommand{\define}[1]{{\bf \boldmath{#1}}}

\theoremstyle{definition}

        \newcommand{\be}{\begin{equation}}
        \newcommand{\ee}{\end{equation}}
        \newcommand{\ba}{\begin{eqnarray}}
        \newcommand{\ea}{\end{eqnarray}}
        \newcommand{\ban}{\begin{eqnarray*}}
        \newcommand{\ean}{\end{eqnarray*}}
        \newcommand{\barr}{\begin{array}}
        \newcommand{\earr}{\end{array}}

\begin{document}
\title{The Gauss--Lucas Theorem}
\author[Baez]{John C.\ Baez} 
\address{Department of Mathematics, University of California, Riverside CA, 92521, USA}
\address{Centre for Quantum Technologies, National University of Singapore, 117543, Singapore}
\date{\today}
\maketitle

The Gauss--Lucas theorem says that for any complex polynomial $P$, the roots of the derivative $P'$ lie in the convex hull of the roots of $P$.   In other words, the roots of $P'$ lie inside the smallest convex subset of the complex plane containing all the roots of $P$.   This theorem is not hard to prove, but is there an intuitive explanation?  In fact there is, using physics---or more precisely, electrostatics in 2-dimensional space \cite{Marden}.

Here is the basic idea.  For simplicity, suppose $P$ has no repeated roots.   Put a particle of charge $1$ at each root.   These particles create a vector field called the electric field, $\E$.   The roots of $P'$ are precisely the points where this electric field vanishes!    In fact, each particle creates its own electric field pointing radially outward, and we sum these fields to get $\E$.   Because of this, $\E$ can only vanish inside the convex hull of the roots of $P$.

\vskip 1em
\begin{center}
\includegraphics[width = 20em]{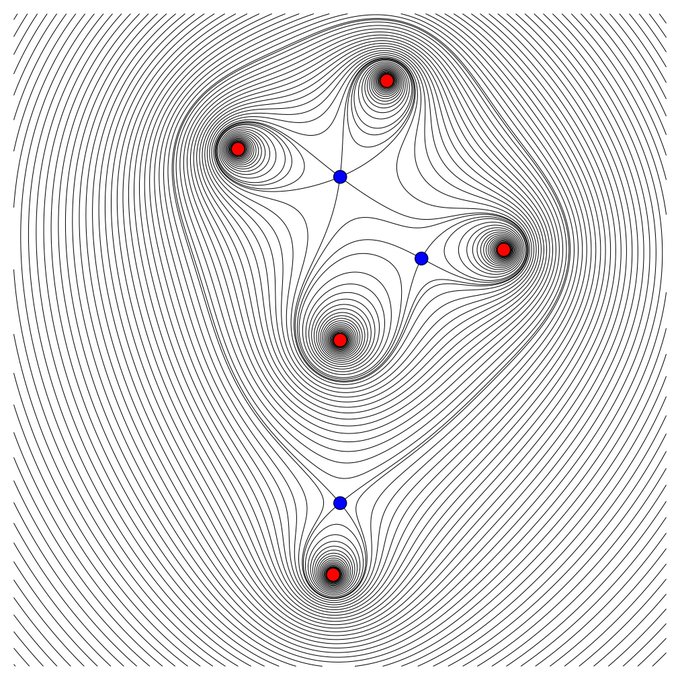} \break
Figure 1 
\end{center}
\vskip 1em

Figure 1, drawn by Greg Egan, shows an example.   Here we see five charged particles marked
in red, located at the roots of some quintic polynomial $P$.  The curves are the level curves of a function $\phi$, called the \define{potential}, such that $\E = -\vec\nabla \phi$.  The potential has three critical points, marked in blue.  The electric field vanishes only at these points.   Note that in some sense these points lie ``between'' the charged particles.  This is in accord with the Gauss--Lucas theorem, and also physical intuition: since each particle creates its own electric field pointing radially outwards, these fields can only cancel at locations between particles.   Note also that our polynomial $P$ has five roots, but its derivative has only three.   But there is no paradox here: in this example, the uppermost blue dot is a root of multiplicity two for $P'$.   Thus, $P'$ has degree four.

To prove the Gauss--Lucas theorem using these ideas, we need some mathematics coming from physics.   The equations of electrostatics make sense in any dimension.  In $\R^n$ they say the potential $\phi \maps \R^n \to \R$ obeys $\nabla^2 \phi = -\rho$ where $\rho \maps \R^n \to \R$ is the charge density.    The case $n = 2$ is special because we can identify $\R^2$ with the complex plane.   In this case, when $\rho$ is a Dirac delta at the origin, there is a solution
\[       \phi(z) = -\frac{1}{2\pi} \ln|z|  .\]
This gives an electric field $\E$ pointing radially outward, with magnitude $1/2\pi|z|$.   This is the electric field produced by a point particle of charge $1$ at the origin.   Since the equations of electrostatics are linear and translation-invariant, it is then easy to solve them for any collection of charged particles.  If we have particles of charge $1$ at points $a_1, \dots, a_k \in \C$, the potential is
\[
 \phi(z) = -\frac{1}{2\pi} \sum_{i = 1}^k \ln |z - a_i|   .
\]
 
Now suppose $P$ is a polynomial without repeated roots, having roots at the points $a_1, \dots, a_k\in \C$.   Putting particles of charge $1$ at these points, we get the above potential $\phi$.  Thus, $|P|$ is some nonzero constant times 
\[    \prod_{i = 1}^k |z - a_i| =  \exp \left(\sum_{i = 1}^k \ln |z - a_i| \right) =  \exp(-2 \pi \phi). \]
The critical points of $P$ are the same as those of $|P|$, except for critical points where $P$  vanishes---but these are forbidden, since we are assuming $P$ has no repeated roots.   So, the critical points of $P$ are precisely the critical points of $\exp(-2 \pi \phi)$.  These, in turn, are the same as the critical points of $\phi$.   But since $\phi = -\nabla \E$, these are the points where the electric field vanishes.   In short, $P'$ vanishes precisely where $\E$ vanishes.

Now let us use this to prove the Gauss--Lucas theorem.   Say $z \in \C$ lies outside the convex hull of the roots $a_1, \dots, a_k$.   Think of all these points as points $\a_i$ and $\z$ in $\R^2$.  By the separating hyperplane theorem, we can draw a line with $\z$ on one side and all the points $\a_i$ on the other side.   Let $\v$ be a vector perpendicular to this line and pointing toward the $\z$ side.   The electric field $\E$ at $\z$ is a sum of vectors pointing from the charges $\a_i$ to $\z$, since
\[    \E(\z) = - \vec\nabla \phi(\z) = 
\frac{1}{2\pi} \sum_{i = 1}^k \frac{\z - \a_i}{|\z - \a_i|^2}  \,.\]
But  $(\z - \a_i) \cdot \v > 0$, so $\E(\z)$ cannot vanish.
  
This proves the Gauss--Lucas theorem for polynomials without repeated roots.  But we only used the assumption that all the $a_i$ were distinct at one place in the argument: to avoid situations where both $P$ and $P'$ vanish at the same point.  These situations pose no problem: a root of $P'$ that is also a root of $P$ is \emph{obviously} in the convex hull of the roots of $P$.

For more thoughts along these lines, we recommend the book by Marden \cite{Marden}.   For example, in the late 1800s, Siebeck, Lucas and B\^{o}cher independently showed that if $P$ is a cubic whose vertices form a triangle, the roots of $P'$ lie at foci of an ellipse inscribed in this triangle---and in fact, this ellipse touches the midpoints of the triangle's sides!   It is called the \define{Steiner inellipse}, and it is also the ellipse of largest area that fits inside the triangle.    
B\^{o}cher asked if this result could be generalized to polynomials of higher degree.  In 1919, Linfield showed that indeed it can, using a general concept of ``foci'' of a real plane algebraic curve \cite{Linfield}.  It would be nice if there were a way to understand some of these more advanced results using electrostatics.

\end{document}